\documentclass[a4paper,10pt]{report}
\usepackage{amssymb}
\usepackage{latexsym}
\usepackage{graphicx}
\usepackage{natbib}
\usepackage{amsmath}

\def\lemma{\textbf{Lemma}}

\def\pf{\textbf{Proof }}
\def\thm{\textbf{Theorem}}

\def\H2{\B{H}^{2}}
\def\diamond{\diamondsuit}

\newcommand{\lra}{\longrightarrow}
\newcommand{\nt}{\newtheorem}
\newcommand{\tit}{\textit}
\newcommand{\tbf}{\textbf}

\newcommand{\B}{\Bbb}

\begin{document}
\begin{center}
\tbf{Combinatorial rigidity in curve complexes and mapping class groups.}

Kenneth J. Shackleton\footnote{Supported by an EPSRC PhD scholarship.

}

e-mail: kjs@maths.soton.ac.uk

http://www.maths.soton.ac.uk/$\sim$kjs

[To appear in the Pacific Journal of Mathematics.]\\
\end{center}

\setlength{\parindent}{0em}

A\begin{small}BSTRACT\end{small}: In all possible cases, we prove that local embeddings between two curve complexes whose complexities do not increase from domain to codomain are induced by surface homeomorphism. This is our first main result. From this we can deduce our second, a strong local co-Hopfian result for mapping class groups.\\

K\begin{small}EYWORDS\end{small}: Curve complex, mapping class group.\\\\

\tbf{$\S0$. Introduction.\\}

The curve complex $\mathcal{C}(\Sigma)$ associated to a surface $\Sigma$ was introduced by Harvey [H] to encode the large scale geometry of Teichm\"uller space, and help decide the non-arithmeticity of the surface mapping class groups. It was later to play a central role in the proof of Brock-Canary-Minsky [BroCMin] of Thurston's ending lamination conjecture.

\setlength{\parindent}{2em}

We start by defining the curve complex, and throughout our surfaces will be compact, connected and orientable. We say that a simple loop on $\Sigma$ is \tit{trivial} if it bounds a disc and \tit{peripheral} if it bounds an annulus whose other boundary component belongs to $\partial \Sigma$. A \tit{curve} on $\Sigma$ is a free homotopy class of a non-trivial and non-peripheral simple loop and we denote the set of these by $X(\Sigma)$. The \tit{intersection number} of two curves $\alpha, \beta \in X(\Sigma)$, denoted $\iota(\alpha, \beta)$, is defined equal to ${\rm min}\{|a \cap b| : a \in \alpha, b \in \beta\}$. We say that two curves \tit{intersect minimally} if they intersect once or they intersect twice with zero algebraic intersection and refer to either as the type of minimal intersection. We will later define the \tit{complexity of $\Sigma$}, denoted $\kappa(\Sigma)$, as equal to the maximal number of distinct and disjoint curves that can be realised simultaneously.

When $\kappa(\Sigma) \geq 2$, the curve graph is the graph whose vertex set is $X(\Sigma)$ and we deem two distinct curves to span an edge if and only if they can be realised disjointly in $\Sigma$. When $\kappa(\Sigma) = 1$, we say that two distinct curves are joined by an edge if and only if they intersect minimally. The \tit{curve complex} associated to $\Sigma$ is the curve graph when $\kappa(\Sigma) = 1$, making it isomorphic to a Farey graph, and the flag simplicial complex whose $1$-skeleton is the curve graph when $\kappa(\Sigma) \geq 2$. In the latter case, $\mathcal{C}(\Sigma)$ has simplicial dimension precisely $\kappa(\Sigma) - 1$. 

For each curve $\alpha$ we denote by $X(\alpha)$ the set of all curves on $\Sigma$ distinct and disjoint from $\alpha$, that is the vertex set of the link of $\alpha$. This link is always connected whenever $\kappa(\Sigma)$ is at least three, and whenever $\kappa(\Sigma)$ is two any two elements of $X(\alpha)$ may be ``chain-connected'' by a finite sequence of curves in which any two consecutive curves have minimal intersection.

In this paper, we shall be discussing embeddings between two curve complexes whose complexities do not increase from domain to codomain and we shall find that these are all induced by surface homeomorphism, so long as we place a necessary but consistent hypothesis in one sporadic case. The argument we give is by an induction on complexity and requires little more than the connectivity of links in the curve complex over and above this. As such, our approach does not discriminate in terms of the topological type of a surface. Moreover, we actually only require the local injectivity of an embedding and we shall say more on this towards the end of this section.

All told, this generalises the automorphism theorem of Ivanov's for surfaces of genus at least two, a proof of which is sketched in $\S2$ of [Iva1] and extended by Korkmaz [K] to all but the two-holed torus, and of Luo's [L], settled or proven in all cases. Making use of their combined result, Margalit [Mar] establishes the analogue for automorphisms of another important surface complex called the pants complex. There are analogues for other surface complexes, see Schmutz Schaller's [Sch] as one example.

Our first result is stated as follows.\\

\nt{BBBB}{\thm}

\begin{BBBB}
Suppose that $\Sigma_{1}$ and $\Sigma_{2}$ are two compact and orientable surfaces of positive complexity such that the complexity of $\Sigma_{1}$ is at least that of $\Sigma_{2}$, and that when they have equal complexity at most three they are homeomorphic or one is the three-holed torus. Then, any simplicial embedding from $\mathcal{C}(\Sigma_{1})$ to $\mathcal{C}(\Sigma_{2})$ (preserving the separating type of each curve when the two surfaces are homeomorphic to the two-holed torus) is induced by a surface homeomorphism.\\
\end{BBBB}

This covers all possibilities. We remind ourselves that there exist isomorphisms between the curve complex of the closed surface of genus two and the six-holed sphere, the two-holed torus and the five-holed sphere and finally the one-holed torus and the four-holed sphere and that there exists an automorphism of the curve complex associated to the two-holed torus that sends a non-separating curve to an outer curve (see [L] for more details). These are examples of embeddings not induced by a surface homeomorphism. Finally, we point out that there exist embeddings on curve complexes with complexity increasing from domain to codomain not induced by a surface embedding: Easy examples are provided by taking some proper subsurface $\Sigma_{1}$ of $\Sigma_{2}$, and modifying the induced embedding on curve complexes by instead taking just one curve on $\Sigma_{1}$ to a curve on $\Sigma_{2}$ outside of $\Sigma_{1}$.

Among other things, Theorem 1 completes one study of a particular class of self-embedding, initiated by Irmak. This class comprises the superinjective maps, and by definition each preserves the non-zero intersection property of a pair of curves. In [Irm1] the author shows that a superinjective self-map is induced by a surface homeomorphism provided the surface is closed and of genus at least three, in [Irm2] this is extended to non-closed surfaces of genus at least three and surfaces of genus two with at least two holes, and in [Irm3] the author extends this to the remaining two types of genus two surface. Following a now standard strategy, set out by Ivanov, this holds consequences for the mapping class groups of the corresponding surfaces.

The \tit{mapping class group} $Map(\Sigma)$ is the group of all self-homeomorphisms of the surface $\Sigma$, up to homotopy. This is sometimes known as the extended mapping class group, for it contains the group of orientation preserving mapping classes as an index two subgroup. Some of its other subgroups, in particular the Johnson kernel and the Torelli group, are of wide interest (see, respectively, Brendle-Margalit [BreMar] and Farb-Ivanov [FIva], and references contained therein).

The mapping class group has a natural simplicial action on the curve complex, determined by first lifting a curve to a representative loop and then taking the free homotopy class of the image under a representative homeomorphism. The kernel of this action, $Ker(\Sigma)$, is almost always trivial; the only exceptions lie in low complexity, where this kernel is isomorphic to $\B{Z}_{2}$ and generated by the hyperelliptic involution when $\Sigma$ is the one-holed torus, the two-holed torus, or the closed surface of genus two or isomorphic to $\B{Z}_{2} \oplus \B{Z}_{2}$ and generated by two hyperelliptic involutions when the four-holed sphere (this is due to Birman [Birm] and Viro [V]). For a detailed account of the mapping class group and its subgroups, see Ivanov [Iva2] as one place to start.

Theorem 1 implies the following strong co-Hopfian result for mapping class groups. Among other things Theorem 2 has some familiar consequences, namely it follows that the commensurator group of a mapping class group is isomorphic to the same mapping class group and that the outer automorphism group of a mapping class group is trivial. Furthermore, it follows that mapping class groups do not admit a faithful action on another curve complex of no greater dimension and that there can be only one faithful action by any such mapping class group on its curve complex, up to conjugation.\\

\nt{virtcom}[BBBB]{\thm}

\begin{virtcom}
Suppose that $\Sigma_{1}$ and $\Sigma_{2}$ are two compact and orientable surfaces such that the complexity of $\Sigma_{1}$ is at least that of $\Sigma_{2}$ and at least two, and that whenever they both have complexity equal to three they are homeomorphic, though not to the closed surface of genus two, or one is the three-holed torus and that when they both have complexity two they are homeomorphic to the five-holed sphere. Suppose that $H$ is a finite index subgroup of the mapping class group $Map(\Sigma_{1})$. Then, every injection of $H$ into $Map(\Sigma_{2})$ is the restriction of an inner automorphism of $Map(\Sigma_{1})$.\\
\end{virtcom}

The existence of such a homomorphism is to imply the two surfaces are equal. Theorem 2 is a generalisation of a result of Ivanov-McCarthy (Theorem 4 from [IvaMcCar]) where the two authors consider injections defined on mapping class groups associated to surfaces of positive genus.

The combined superinjectivity theorem implies Theorem 2 when the two surfaces under consideration are homeomorphic and have genus at least three, or genus at least two and one hole, and in [Irm3] the author describes a non-inner automorphism for the closed surface of genus two. Bell-Margalit [BelMar2] extend this to spheres with at least five holes, and Behrstock-Margalit [BehrMar] to genus one surfaces with at least three holes in addition to finding a commensurator for the mapping class group of the two-holed torus not induced by an inner automorphism. The remaining cases, namely the mapping class group of the four-holed sphere and of the one-holed torus, also have non-inner injections on finite index subgroups, as both are virtually free groups. We remark the braid groups on at least four strands, modulo centre, are shown by Bell-Margalit [BelMar1] to have the co-Hopfian property.

The general approach we need for Theorem 2 follows that given by Ivanov, translating an injection on a finite index subgroup to an embedding on curve complexes. This is now a well-established strategy on which we have nothing to add, and a thorough account can be found in the work of Bell-Margalit [BelMar2] or [Irm1].


Though all our arguments are phrased in terms of embeddings, they only ever need the simplicial and local injectivity properties of such maps. We can therefore record the following generalisation of Theorem 1, the first of two main results. Recall that a \tit{star} is the union of all edges incident on a common vertex.\\

\nt{folding}[BBBB]{\thm}

\begin{folding}
Suppose that $\Sigma_{1}$ and $\Sigma_{2}$ are two compact and orientable surfaces of positive complexity such that the complexity of $\Sigma_{1}$ is at least that of $\Sigma_{2}$, and that when they have equal complexity at most three they are homeomorphic or one is the three-holed torus. Then, any simplicial map from $\mathcal{C}(\Sigma_{1})$ to $\mathcal{C}(\Sigma_{2})$ injective on every star (and preserving the separating type of each curve when both surfaces are homeomorphic to the two-holed torus) is induced by a surface homeomorphism.\\
\end{folding}

Again, this covers all possibilities. We remark that proving a local embedding is induced by a surface homeomorphism would appear the most direct way of seeing that it must also be a global embedding. Furthermore, we conjecture that the pants complex also exhibits such local-to-global rigidity.

From Theorem 3 we can deduce, using a careful application of Ivanov's strategy, the following local version of Theorem 2. This is one interpretation of local injectivity for mapping class groups, and a proof is completed in $\S3.3$ of [Sha]. Among other things, it follows that a self-homomorphism of a mapping class group injective on every curve stabiliser is the restriction of an inner automorphism.\\

\nt{foldingMCG}[BBBB]{\thm}

\begin{foldingMCG}
Suppose that $\Sigma_{1}$ and $\Sigma_{2}$ are two compact and orientable surfaces such that the complexity of $\Sigma_{1}$ is at least that of $\Sigma_{2}$ and at least three, and that whenever they both have complexity equal to three they are homeomorphic, though not to the closed surface of genus two, or one is the three-holed torus and that whenever they both have complexity two they are homeomorphic to the five-holed sphere. Suppose that $H$ is a finite index subgroup of the mapping class group $Map(\Sigma_{1})$. Then, every homomorphism of $H$ into $Map(\Sigma_{2})$ injective on every curve stabiliser in $H$ is the restriction of an inner automorphism of $Map(\Sigma_{1})$.\\
\end{foldingMCG}

Note once more, the existence of such a homomorphism is to imply the two surfaces are equal.

Investigations into arbitrary homomorphisms from a mapping class group associated to a closed surface of genus at least one to another mapping class group associated to a closed surface of smaller genus have been made by Harvey-Korkmaz [HK], the authors finding that every such homomorphism has finite image. Their approach seems to make essential use of the existence of torsion in mapping class groups and, as mapping class groups are virtually torsion free, it would be of some interest to find a way around this so as to consider finite index subgroups.\\

\setlength{\parindent}{0em}

\tbf{Acknowledgements.} The author wishes to thank the Department of Mathematical and Computing Sciences at the Tokyo Institute of Technology for its warm and generous hospitality, and it is his pleasure to thank Brian Bowditch, Javier Aramayona, Jason Behrstock, Dan Margalit, Graham Niblo and Caroline Series for their enthusiasm and for being so generous with their time; their comments on this work improved the exposition beyond his own capabilities. Last, but by no means least, he wishes to thank the EPSRC for supporting this research.\\\\

\tbf{$\S1$. Embeddings between curve complexes.}\\

For any compact, connected and orientable surface $\Sigma$ the \tit{complexity} $\kappa(\Sigma)$ of $\Sigma$ is defined to be equal to $3genus(\Sigma) + |\partial \Sigma| - 3$. This is perhaps non-standard, since complexity is often taken to be equal to the simplicial dimension of the curve complex, but the additivity of $\kappa$ best suits our induction argument. By way of example, the one-holed torus and the four-holed sphere are the only surfaces of complexity one, the two-holed torus and the five-holed sphere are the only surfaces of complexity two, and the closed surface of genus two, the three-holed torus and the six-holed sphere are the only surfaces of complexity three. On occasion we refer to these as the \tit{low complexity surfaces}.

\setlength{\parindent}{2em}

In what follows, we shall abuse notation slightly by viewing each curve as a vertex, as a class of loops, or as a simple loop already realised on $\Sigma$. Our interpretation will be apparent from the context. We say that a curve is \tit{separating} if its complement is not connected, and otherwise say it is \tit{non-separating}. We say that a curve is an \tit{outer curve} if it is separating and if it bounds a two-holed disc (equivalently, a three-holed sphere). These are usually known as boundary curves in the literature, but here we need to avoid confusing them with the components of $\partial \Sigma$. A \tit{multicurve on $\Sigma$} is a collection of distinct and disjoint curves, and a \tit{pants decomposition of $\Sigma$} is a maximal multicurve. A \tit{pair of pants in $\Sigma$} is an essential subsurface homeomorphic to a compact three-holed sphere.

We say that two curves in a pants decomposition $P$ are \tit{adjacent in $P$} if they appear in the boundary of a non-compact three-holed sphere complementary to $P$. This is slightly unfortunate terminology that seems to be a long way to becoming standard; we sincerely hope that any confusion between adjacency in the curve complex and adjacency in a pants decomposition will be obviated by the context.

The structure of our argument is broadly as follows. We establish a short list of topological properties verified by any embedding on curve complexes from which we easily deduce, among other things, that the existence of such an embedding implies the two surfaces have equal complexity and then, with more work, almost always means the two surfaces under consideration are homeomorphic. For the time being, we refer to embeddings between two apparently distinct curve complexes as \tit{cross-embeddings}. Dealing with embeddings in low complexity typically requires individual arguments and it therefore streamlines our work if we do this separately, as we do in Lemma 13. The proof of Theorem 1 is then completed by an induction on complexity, where we cut the surface along a curve. As embeddings behave well on the topological type of a curve, the resulting surfaces are again homeomorphic. For the induction argument to pass through complexity one (sub)surfaces, we will need to show that embeddings preserve minimal intersection.

We start by showing, in turn, that embeddings send pants decompositions to pants decompositions, they preserve a form of small intersection and they preserve adjacency \tit{and} non-adjacency in a pants decomposition.\\

\nt{panttopant}[BBBB]{\lemma}

\begin{panttopant}
Suppose that $\Sigma_{1}$ and $\Sigma_{2}$ are two compact and orientable surfaces such that the complexity of $\Sigma_{1}$ is at least that of $\Sigma_{2}$. Then, any embedding $\phi$ from $\mathcal{C}(\Sigma_{1})$ to $\mathcal{C}(\Sigma_{2})$ sends pants decompositions to pants decompositions.\\
\end{panttopant}

\setlength{\parindent}{0em}

\pf This follows for complexity reasons and because $\phi$ is simplicial and injective. $\diamond$\\

\setlength{\parindent}{2em}

To make sense of the following lemma, we must define what we mean by the \tit{subsurface of $\Sigma$ filled by two curves $\alpha$ and $\beta$}. Letting $N(\alpha \cup \beta)$ denote a regular neighbourhood of $\alpha \cup \beta$ in $\Sigma$, we augment $N(\alpha \cup \beta)$ by taking its union with all the complementary discs whose boundary is contained in $N(\alpha \cup \beta)$ and all the complementary annuli with one boundary component in $\partial \Sigma$ and the other in $N(\alpha \cup \beta)$. The resulting subsurface of $\Sigma$ is well defined up to homotopy and is what we mean by the subsurface filled by $\alpha$ and $\beta$. Whenever a third curve enters the subsurface filled by two curves, it must intersect at least one of these two curves.\\

\nt{weaksuper}[BBBB]{\lemma}

\begin{weaksuper}
Suppose that $\Sigma_{1}$ and $\Sigma_{2}$ are two compact and orientable surfaces such that the complexity of $\Sigma_{1}$ is at least that of $\Sigma_{2}$. Let $\phi$ be any embedding from $\mathcal{C}(\Sigma_{1})$ to $\mathcal{C}(\Sigma_{2})$ and let $\alpha, \beta$ be any two curves in $\Sigma_{1}$ that fill either a four-holed sphere or a one-holed torus. Then, $\phi(\alpha)$ and $\phi(\beta)$ fill either a four-holed sphere or a one-holed torus in $\Sigma_{2}$.\\
\end{weaksuper}

\setlength{\parindent}{0em}

\pf Let $Q$ be any maximal multicurve in $\Sigma_{1}$ such that each curve is disjoint from both $\alpha$ and $\beta$. For complexity reasons, $\phi(Q)$ is a maximal multicurve disjoint from both $\phi(\alpha)$ and $\phi(\beta)$. In particular, as $\phi$ is injective and simplicial so $\phi(\alpha)$ and $\phi(\beta)$ must together fill either a four-holed sphere or a one-holed torus. $\diamond$\\

\setlength{\parindent}{2em}

We shall say that two curves have \tit{small intersection} if they together fill either a four-holed sphere or a one-holed torus, and refer to either as the \tit{type} of the small intersection. Any two curves that intersect minimally also have small intersection, but the converse does not hold.\\

\nt{presadj}[BBBB]{\lemma}

\begin{presadj}
Suppose that $\Sigma_{1}$ and $\Sigma_{2}$ are two compact and orientable surfaces such that the complexity of $\Sigma_{1}$ is at least that of $\Sigma_{2}$. Let $P$ be any pants decomposition of $\Sigma_{1}$ and let $\phi$ be any embedding from $\mathcal{C}(\Sigma_{1})$ to $\mathcal{C}(\Sigma_{2})$. Then, any two curves adjacent in $P$ are sent by $\phi$ to two curves adjacent in $\phi(P)$ and any two curves in $P$ that are not adjacent in $P$ are sent by $\phi$ to two curves not adjacent in $\phi(P)$.\\
\end{presadj}

\setlength{\parindent}{0em}

\pf The first part follows from Lemma 6: For any two curves $\alpha_{1}$ and $\alpha_{2}$ adjacent in $P$, there exists a curve $\delta$ having small intersection with both and disjoint from every other curve in $P$. This is preserved under $\phi$ and so $\phi(\alpha_{1})$ and $\phi(\alpha_{2})$ are adjacent in $\phi(P)$. 

\setlength{\parindent}{2em}

Similarly, if two curves $\alpha_{1}, \alpha_{2}$ are not adjacent in $P$ we can find two disjoint curves $\delta_{1}, \delta_{2}$ such that $\delta_{1}$ has small intersection with $\alpha_{1}$ but is disjoint from $\alpha_{2}$ and $\delta_{2}$ has small intersection with $\alpha_{2}$ but is disjoint from $\alpha_{1}$ and both $\delta_{1}, \delta_{2}$ are disjoint from every other curve in $P$. If $\phi(\alpha_{1})$ and $\phi(\alpha_{2})$ are adjacent in $\phi(P)$ then $\phi(\delta_{1})$ and $\phi(\delta_{2})$ must intersect. As $\phi$ is simplicial, this is a contradiction. $\diamond$\\

The import of Lemma 5, Lemma 6 and Lemma 7 is perhaps best understood by associating to a pants decomposition $P$ a certain graph. The vertices of this graph are the curves in $P$, and any two distinct vertices span an edge if and only if they correspond to adjacent curves in $P$. Lemma 7 not only tells us that any embedding $\phi$ induces a map between adjacency graphs, but that this map is actually an isomorphism. Cut points in the graph correspond to non-outer separating curves, and non-cut points correspond to outer or non-separating curves.

This graph, and the ideas bound by Lemma 7, were independently and simultaneously discovered by Behrstock-Margalit. Their approach can be found in [BehrMar] and the arguments they give will deal with all superinjective maps for two homeomorphic surfaces of complexity at least three. From this they also deduce that the commensurator group of a mapping class group is isomorphic to the same mapping class group. We both refer to such a graph as an \tit{adjacency graph}.

\begin{figure}
\centering 
\includegraphics[width=0.75\textwidth]{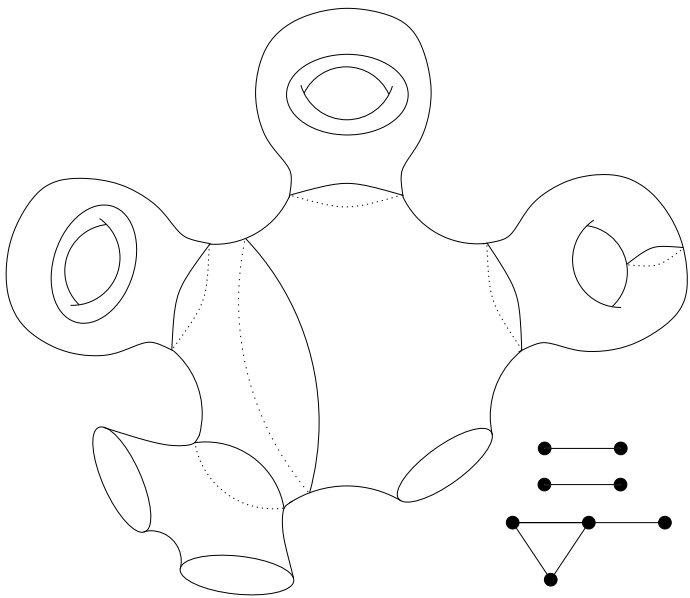} 
\caption{A codimension $1$ multicurve, with its adjacency graph.} 
\label{adj1}
\end{figure}

We can just as well speak of an adjacency graph associated to a multicurve $Q$, in which the vertices again correspond to the curves in $Q$ and any two vertices are declared adjacent if their corresponding curves border a common pair of pants in the surface complement of $Q$. There is a subtle point to be made here, namely that the complementary graph of a vertex in a pants adjacency graph will \tit{not} in general be the adjacency graph of the multicurve that results by removing the corresponding curve from the pants decomposition. It will however be the adjacency graph that results from cutting the surface along this curve. By way of example, on removing a curve $\alpha$ from a pants decomposition $P$ any curves that together bound a complementary four-holed sphere will not necessarily be adjacent in the adjacency graph of $P-\{\alpha\}$. (See Figure \ref{adj1} for an illustrated example.) This observation will be important later when we come to look at outer curves. It does however hold that a curve complex embedding induces an isomorphism between multicurve adjacency graphs.\\

\nt{adjmulti}[BBBB]{\lemma}

\begin{adjmulti}
Suppose that $\Sigma_{1}$ and $\Sigma_{2}$ are two compact and orientable surfaces such that the complexity of $\Sigma_{1}$ is at least that of $\Sigma_{2}$. Let $Q$ be any multicurve of $\Sigma_{1}$ and let $\phi$ be any embedding from $\mathcal{C}(\Sigma_{1})$ to $\mathcal{C}(\Sigma_{2})$. Then, $\phi$ induces an isomorphism from the adjacency graph of $Q$ to the adjacency graph of $\phi(Q)$.\\
\end{adjmulti}

\setlength{\parindent}{0em}

\pf We make use of Lemma 7. Extend $Q$ to a pants decomposition $P$ of $\Sigma_{1}$. If two curves are adjacent in $Q$ then they either border a pair of pants with a third curve from $Q$ or they border a pair of pants meeting $\partial \Sigma$. This remains so in $P$, and is preserved on applying $\phi$. To show non-adjacency is preserved, consider any two curves not adjacent in $Q$ and arrange for them to be non-adjacent in $P$. This is preserved under $\phi$. $\diamond$\\

\setlength{\parindent}{2em}

As embeddings between curve complexes induce isomorphisms on adjacency graphs and graph isomorphisms send cut points to cut points, so embeddings must send non-outer separating curves to non-outer separating curves.\\

\nt{septosep}[BBBB]{\lemma}

\begin{septosep}
Suppose that $\Sigma_{1}$ and $\Sigma_{2}$ are two compact and orientable surfaces such that the complexity of $\Sigma_{1}$ is at least that of $\Sigma_{2}$. Then, any embedding $\phi$ from $\mathcal{C}(\Sigma_{1})$ to $\mathcal{C}(\Sigma_{2})$ sends non-outer separating curves to non-outer separating curves.\\
\end{septosep}

We use Lemma 6 and the adjacency graph to distinguish between non-separating and outer curves.\\

\nt{nontonon}[BBBB]{\lemma}

\begin{nontonon}
Suppose $\Sigma_{1}$ and $\Sigma_{2}$ are two compact and orientable surfaces such that the complexity of $\Sigma_{1}$ is at least that of $\Sigma_{2}$ and that whenever they have equal complexity at most three they are homeomorphic and not the two-holed torus. Let $\phi : \mathcal{C}(\Sigma_{1}) \lra \mathcal{C}(\Sigma_{2})$ be any embedding. Then, $\phi$ takes non-separating curves to non-separating curves.\\
\end{nontonon}

\setlength{\parindent}{0em}

\pf We note that the $\phi$-image of a non-separating curve can never be a non-outer separating curve, for otherwise we see a non-cut point sent to a cut point in some pants adjacency graph. Suppose that $\alpha$ is a non-separating curve in $\Sigma_{1}$. When $\kappa(\Sigma_{1})$ is at least four we can find a pants decomposition $P$ extending $\alpha$ in which $\alpha$ corresponds to a vertex in the adjacency graph of $P$ of valence three or four. As $\phi$ induces an isomorphism on the adjacency graph, so $\phi(\alpha)$ must have the same valence. As outer curves only ever correspond to vertices of valence at most two, so $\phi(\alpha)$ can only be non-separating.

\setlength{\parindent}{2em}

With the exception of the two-holed torus, all cases in which $\Sigma_{1}$ has complexity at most two hold since there is only ever one type of curve. In complexity three, when $\Sigma_{1}$ is the six-holed sphere our claim holds vacuously and when $\Sigma_{1}$ is the closed surface of genus two our claim follows from Lemma 9 by noting that every pants decomposition contains at most one separating curve. 

The only non-trivial case in low complexity is that of $\Sigma_{1}$ and $\Sigma_{2}$ both homeomorphic to the three-holed torus. In which case, there are only two pants adjacency graphs, up to isomorphism, but three different pants decompositions, up to the action of the mapping class group. For this reason, we need to argue differently. If there is a non-separating curve sent by $\phi$ to an outer curve, then there is an outer curve $\alpha$ sent by $\phi$ to a non-separating curve. To see this, extend this non-separating curve to a pants decomposition containing a non-outer separating curve. By appealing to Lemma 9, we see that the third curve in this pants decomposition will suffice. Now extend $\alpha$ to a second pants decomposition containing two non-separating curves $\delta_{1}$ and $\delta_{2}$. The $\phi$-image of at least one of these, say $\delta_{1}$, is again a non-separating curve. Choose any two disjoint curves $\gamma_{1}, \gamma_{2}$ in $\Sigma_{1}$ that have small intersection with $\delta_{1}$ and $\alpha$ but disjoint from $\alpha$ and $\delta_{1}$, respectively. Now $\phi(\delta_{1})$ and $\phi(\alpha)$ border a common pair of pants in $\Sigma_{2}$ invaded by $\phi(\gamma_{1})$ and $\phi(\gamma_{2})$. We see that the $\phi$-images of both $\gamma_{1}$ and $\gamma_{2}$ are forced to intersect, and this is a contradiction. $\diamond$\\

\nt{boundtobound}[BBBB]{\lemma}

\begin{boundtobound}
Suppose that $\Sigma_{1}$ and $\Sigma_{2}$ are two compact and orientable surfaces such that the complexity of $\Sigma_{1}$ is at least that of $\Sigma_{2}$, and that whenever they have equal complexity at most three they are homeomorphic and not the two-holed torus. Then, any embedding $\phi$ from $\mathcal{C}(\Sigma_{1})$ to $\mathcal{C}(\Sigma_{2})$ sends outer curves to outer curves.\\
\end{boundtobound}

\setlength{\parindent}{0em}

\begin{figure}
\centering 
\includegraphics[width=0.75\textwidth]{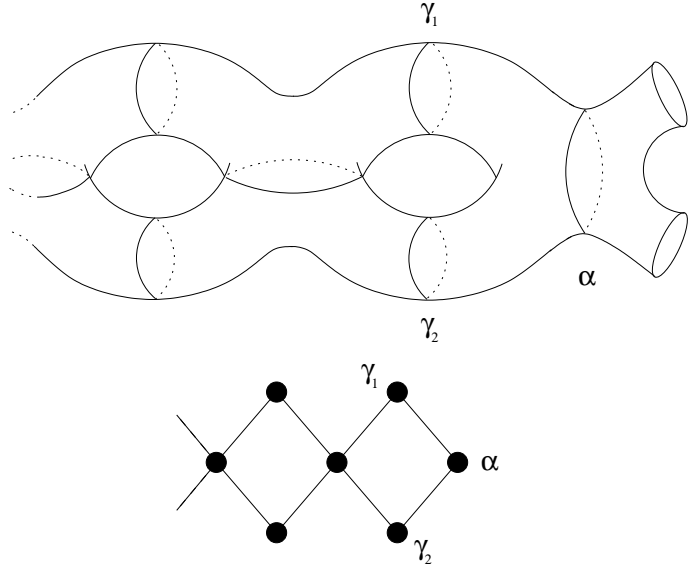} 
\caption[A convenient choice of pants decomposition extending an outer curve.]{A convenient extension of $\alpha$ to a pants decomposition.} 
\label{choice}
\end{figure}

\pf We note that this holds vacuously when $|\partial \Sigma_{1}|$ is at most one. In any case, let us suppose for contradiction that $\alpha$ is an outer curve in $\Sigma_{1}$ sent by $\phi$ to a non-outer curve in $\Sigma_{2}$. We note that $\phi(\alpha)$ can not be a separating curve, for $\alpha$ can never correspond to a cut point in a pants adjacency graph, and so $\phi(\alpha)$ must be a non-separating curve. If $\kappa(\Sigma_{1})$ is at least four then we can extend $\alpha$ to a pants decomposition $P$ in which the two curves adjacent to $\alpha$, denoted $\gamma_{1}$ and $\gamma_{2}$, are not adjacent in the adjacency graph of $P - \{\alpha\}$. As $\alpha$ is an outer curve, we note that $\gamma_{1}$ and $\gamma_{2}$ are adjacent in $P$. According to Lemma 8, $\phi(\gamma_{1})$ and $\phi(\gamma_{2})$ can only, together with $\partial \Sigma_{2}$, border a four-holed sphere containing $\phi(\alpha)$. However, $\phi(\alpha)$ is not an outer curve. Therefore, $\phi(\gamma_{1})$ and $\phi(\gamma_{2})$ are not adjacent in $\phi(P)$ and this is contrary to the statement of Lemma 7. (See Figure \ref{choice} for one example.)

\setlength{\parindent}{2em}

Once more, the only remaining non-trivial case in low complexity is that of the three-holed torus. Suppose that $\alpha$ is an outer curve sent to a non-separating curve by $\phi$. Extend $\alpha$ to a pants decomposition $P$ containing a separating curve. Then the non-separating curve in $P$ is sent to an outer curve by $\phi$, and this is contrary to Lemma 10. $\diamond$\\

It now follows that small subsurfaces can not change topological type under embeddings.\\

\nt{nochange}[BBBB]{\lemma}

\begin{nochange}
Suppose that $\Sigma_{1}$ and $\Sigma_{2}$ are two compact and orientable surfaces such that the complexity of $\Sigma_{1}$ is at least that of $\Sigma_{2}$ and that when both have equal complexity at most three they are homeomorphic and not the two-holed torus. Let $Z$ be any essential $\pi_{1}$-injective subsurface of $\Sigma_{1}$ of complexity one and bordered by a single curve $\beta$. Then, for any embedding $\phi$ from $\mathcal{C}(\Sigma_{1})$ to $\mathcal{C}(\Sigma_{2})$, the $\pi_{1}$-injective minimal subsurface $\phi(Z)$ of $\Sigma_{2}$ filled by $\phi(X(Z))$ is homeomorphic to $Z$.\\
\end{nochange}

\setlength{\parindent}{0em}

\pf Such a change in topology would otherwise force $\phi$ to send a non-separating curve to an outer curve or an outer curve to a non-separating, contrary to Lemma 10 and Lemma 11 respectively. $\diamond$\\

\setlength{\parindent}{2em}

We can finally rule out cross-embeddings, and thereafter we regard the two surfaces as being equal and denote both by $\Sigma$.\\ 

\nt{nolowcross}[BBBB]{\lemma}

\begin{nolowcross}
Suppose that $\Sigma_{1}$ and $\Sigma_{2}$ are two compact and orientable surfaces such that the complexity of $\Sigma_{1}$ is at least that of $\Sigma_{2}$, and that whenever they have complexity at most two they are homeomorphic and whenever they have complexity equal to three they are either homeomorphic or one is the three-holed torus. Then, there is an embedding $\phi : \mathcal{C}(\Sigma_{1}) \lra \mathcal{C}(\Sigma_{2})$ only if $\Sigma_{1}$ and $\Sigma_{2}$ are homeomorphic.\\
\end{nolowcross}

\setlength{\parindent}{0em}

\pf The existence of such an embedding implies the complexities $\kappa(\Sigma_{1})$ and $\kappa(\Sigma_{2})$ are equal. When $\kappa(\Sigma_{1})$ is at least four, we know that any such embedding must send separating curves to separating curves. We recall that the size of a maximal collection of distinct and disjoint separating curves in $\Sigma_{1}$ is precisely $2genus(\Sigma_{1}) + |\partial \Sigma_{1}| - 3$. By our earlier work, this is at most $2genus(\Sigma_{2}) + |\partial \Sigma_{2}| - 3$ and so $genus(\Sigma_{1}) \geq genus(\Sigma_{2})$.

\setlength{\parindent}{2em}

To prove equality, we take $Q$ to be a maximal collection of distinct and disjoint curves on $\Sigma_{1}$ each bounding a one-holed torus. That is, $Q$ has $genus(\Sigma_{1})$ curves. According to Lemma 12, each curve in $\phi(Q)$ must also bound a one-holed torus in $\Sigma_{2}$. We deduce $genus(\Sigma_{1}) \leq genus(\Sigma_{2})$. Combining the two inequalities we have $genus(\Sigma_{1}) = genus(\Sigma_{2})$, and $\Sigma_{1}$ and $\Sigma_{2}$ are homeomorphic.

Turning to the low complexity surfaces, there are no embeddings from the curve complex of the six-holed sphere or closed surface of genus two to that of the three-holed torus. To see this, extend an outer or non-separating curve $\alpha$ in $\Sigma_{1}$ to a pants decomposition $P$ consisting only of outer or non-separating curves, respectively, and choose a separating curve $\beta$ disjoint from both curves in $P - \{\alpha\}$ and therefore of small intersection with $\alpha$. We may assume that if any curve in $\phi(P)$ is outer then it is $\phi(\alpha)$. Now $\phi(\beta)$ is a non-outer separating curve intersecting $\phi(\alpha)$ and it follows that $\phi(\beta)$ must intersect another curve in $\phi(P)$. This is a contradiction. 

The remaining cases, namely from the curve complex of the three-holed torus to the curve complex of the six-holed sphere or of the closed surface of genus two, are covered as follows: For any pants decomposition $P$ in $\Sigma_{1}$ comprising only of non-separating curves, choose a non-outer separating curve $\beta$ meeting only two curves in $P$. When $\Sigma_{2}$ is the six-holed sphere, according to Lemma 8 each curve in $P$ can only go to an outer curve. By Lemma 6, small intersection is preserved. Now any non-outer separating curve in the six-holed sphere meets either only one curve or all three curves in a pants decomposition made up entirely of outer curves. It follows that $\phi(\beta)$ meets every curve in $\phi(P)$, and this is a contradiction. This simultaneously deals with $\Sigma_{2}$ the closed surface of genus two. $\diamond$\\

To allow the induction argument to pass through complexity one surfaces unhindered, we need the following lemma on minimal intersection in those subsurfaces bordered by a single curve. This relies on what is a well-established argument, given by Ivanov [Iva1] for intersection one and by Luo [L] for intersection two with zero algebraic intersection. Although both are stated for automorphisms, both apply in our setting.\\

\nt{pres1or20}[BBBB]{\lemma}

\begin{pres1or20}
Suppose that $\Sigma$ is a compact and orientable surface of positive complexity and not homeomorphic to the two-holed torus. Suppose that $Z$ is an essential subsurface of $\Sigma$ of complexity one and bordered by a single curve $\beta$. Then, any embedding $\phi : \mathcal{C}(\Sigma) \lra \mathcal{C}(\Sigma)$ preserves minimal intersection and its type on $X(Z)$.\\
\end{pres1or20}



\setlength{\parindent}{2em}


This closes our study of the topological properties of curve complex embeddings, and the promised induction argument now starts with a look at the Farey graph.\\

\nt{farey}[BBBB]{\lemma}

\begin{farey}
Every simplicial embedding from a Farey graph $\mathcal{F}$ to itself is an automorphism.\\
\end{farey}

\setlength{\parindent}{0em}

\pf We note that each edge in $\mathcal{F}$ separates and belongs to exactly two $3$-cycles and that such a map sends $3$-cycles to $3$-cycles. Thus, any embedding $\phi$ on $\mathcal{F}$ induces an embedding $\phi^{*}$ on the dual graph. This graph is a tree in which every vertex has the same valence, hence the induced map is a surjection. It follows that every $3$-cycle of $\mathcal{F}$ is contained in the image of $\phi$. That is to say, $\phi$ is also a surjection. $\diamond$\\

\setlength{\parindent}{2em}

It is a well-known fact (indeed, it was known to Dehn [D]) that the automorphisms of $\mathcal{C}(\Sigma)$ are all induced by surface homeomorphisms when $\Sigma$ is either a four-holed sphere or a one-holed torus. This completes the base case of the induction.

We now furnish the inductive step. Let $\phi : \mathcal{C}(\Sigma) \lra \mathcal{C}(\Sigma)$ be any embedding satisfying the hypotheses of Theorem 1. Let $\alpha$ be any curve in $\Sigma$. Our previous work on the topological properties of $\phi$ tells us that the complement of $\alpha$ and the complement of $\phi(\alpha)$ are homeomorphic. Therefore, after first composing with a suitable mapping class if need be, $\phi$ restricts to a self-embedding on the curve complex associated to each component of $\Sigma-\alpha$. The embeddings arising in this way are very natural for they inherit many of the properties verified by $\phi$, for instance they also preserve the separating type of a curve. This is of particular relevance when cutting the surface $\Sigma$ along a curve and finding a two-holed torus complementary component. In [L], the author explains how to find automorphisms of the curve complex associated to the two-holed torus not induced by a surface homeomorphism. No such automorphism can arise as a restriction, nor can any embedding, as outer curves in this two-holed torus correspond to separating curves in $\Sigma$. 

Our inductive hypothesis therefore applies and it tells us that each restriction of $\phi$ associated to a positive complexity component of $\Sigma - \alpha$ is induced by a surface homeomorphism. In gluing back together by identifying the boundary components of $\Sigma - \alpha$ corresponding to $\alpha$, we have a countable family of mapping classes where each such mapping class $f$ satisfies $f(\delta) = \phi(\delta)$ for all $\delta \in X(\alpha) \cup \{\alpha\}$. We must somehow decide which of these, if any, is appropriate.

This construction applies equally well for every curve on $\Sigma$, in particular any curve $\beta$ adjacent to $\alpha$. The set of mapping classes associated to $\alpha$ and the set of mapping classes associated to $\beta$ have non-empty intersection. That is, to the edge of $\mathcal{C}(\Sigma)$ spanned by $\alpha$ and $\beta$, we can associate at least one mapping class $f$ with $f(\delta) = \phi(\delta)$ for all $\delta \in X(\alpha) \cup X(\beta)$.

We need to verify that for any three curves $\alpha, \beta_{1}$ and $\beta_{2}$ such that $\alpha$ is adjacent to both $\beta_{1}$ and $\beta_{2}$, the action on $\mathcal{C}(\Sigma)$ of any such mapping class $f_{1}$ associated to the edge $\alpha, \beta_{1}$ agrees with that of any such mapping class $f_{2}$ associated to the edge $\alpha, \beta_{2}$. For almost all surfaces, $f_{1}$ and $f_{2}$ will be the same mapping class. For now there remains the possibility that $\Sigma - \alpha$ has an exceptional surface component and the possibility that $f_{1}^{-1}f_{2}$ Dehn twists around $\alpha$, or a combination of the two. We treat this in the following lemma.\\

\nt{glue}[BBBB]{\lemma}

\begin{glue}
Suppose that $\alpha, \beta_{1}, \beta_{2} \in X(\Sigma)$ are distinct, with $\beta_{1}$ and $\beta_{2}$ of zero or otherwise minimal intersection and $\alpha$ disjoint from both $\beta_{1}$ and $\beta_{2}$. Suppose $f_{1}, f_{2} \in Map(\Sigma)$ are two mapping classes such that $f_{i}(\delta) = \phi(\delta)$ for all $\delta \in X(\alpha) \cup X(\beta_{i})$, for $i = 1, 2$. Then, $f_{1}^{-1}f_{2} \in Ker(\Sigma)$.\\
\end{glue}

\setlength{\parindent}{0em}

\pf Let $f$ denote the mapping class $f_{1}^{-1}f_{2}$, noting $f$ acts trivially on $X(\alpha)$, and suppose for contradiction that $f \notin Ker(\Sigma)$. As we shall see in the subsequent paragraphs, there then exist disjoint (possibly equal) curves $\delta_{1}$ and $\delta_{2}$ on $\Sigma$ such that at least one of $\iota(\delta_{i}, \alpha)$ and $\iota(\delta_{i}, \beta_{i})$ is zero, for both $i = 1, 2$, and such that $\iota(\delta_{1}, f(\delta_{2})) > 0$. Given this, we also have $\iota(\delta_{1}, f(\delta_{2})) = \iota(\delta_{1}, f_{1}^{-1}f_{2}(\delta_{2})) = \iota(f_{1}(\delta_{1}), f_{2}(\delta_{2})) = \iota(\phi(\delta_{1}), \phi(\delta_{2})) = 0$. This is a contradiction, and we deduce the statement of the lemma.

\setlength{\parindent}{2em}

\begin{figure}
\centering 
\includegraphics[width=0.75\textwidth]{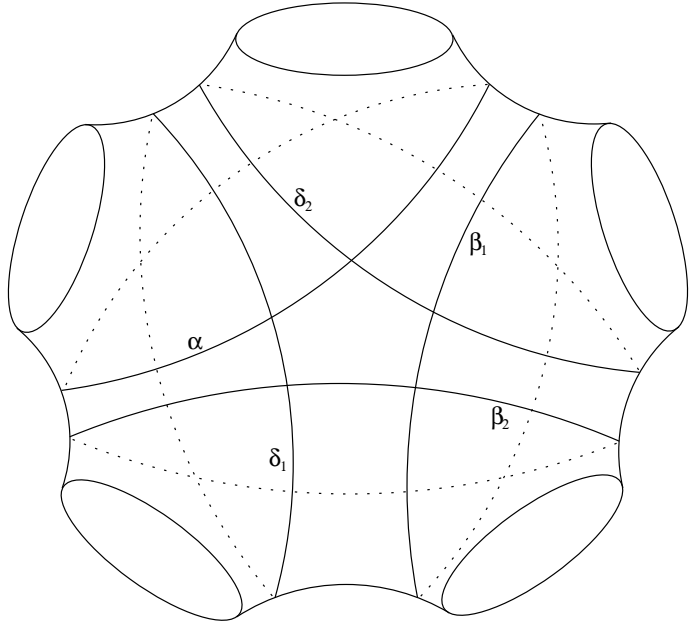} 
\caption{The case $\Sigma$ a five-holed sphere. Note $\delta_{1}$ is the only curve on $\Sigma$ disjoint from both $\beta_{1}, \delta_{2} \in X(\delta_{1}) \cap (X(\alpha) \cup X(\beta_{2}))$. As it happens, $X(\delta_{1}) \cap (X(\alpha) \cup X(\beta_{2}))$ is precisely $\{\beta_{1}, \delta_{2}\}$ in this instance.} 
\label{pent}
\end{figure}

To see that such a pair of curves $\delta_{1}, \delta_{2}$ must exist, we can argue as follows. Suppose $\delta_{1} \in X(\beta_{1})$ has minimal intersection with $\alpha$ and zero or minimal intersection with $\beta_{2}$. Then, $\delta_{1}$ is entirely determined by the non-empty set $X(\delta_{1}) \cap (X(\alpha) \cup X(\beta_{2}))$. More precisely, $\delta_{1}$ is the only curve on $\Sigma$ intersecting $\alpha$ and disjoint from every curve in $X(\delta_{1}) \cap (X(\alpha) \cup X(\beta_{2}))$. See Figure \ref{pent} for one example in the five-holed sphere, in this case a pentagon configuration as described in [L].


Suppose for contradiction that $\iota(\delta_{1}, f(\delta_{2})) = 0$ for any curve $\delta_{2} \in X(\delta_{1}) \cap (X(\alpha) \cup X(\beta_{2}))$. Then, $f(X(\delta_{1}) \cap (X(\alpha) \cup X(\beta_{2}))) \subseteq X(\delta_{1}) \cap (X(\alpha) \cup X(\beta_{2}))$. However, because $f$ is a mapping class this inclusion is an equality and we deduce $f(\delta_{1}) = \delta_{1}$. As the complement in $\Sigma$ of $\beta_{1}$ is filled by a set of curves all fixed by $f$, we deduce $f$ acts trivially on $X(\beta_{1})$. Arguing along similar lines, by reinterpreting our contention as $\iota(\delta_{2}, f^{-1}(\delta_{1})) = 0$ we deduce $f$ acts trivially on $X(\beta_{2})$ as well.

We have shown that $f(\delta) = \delta$ for all $\delta \in X(\alpha) \cup X(\beta_{1}) \cup X(\beta_{2})$. However $X(\alpha) \cup X(\beta_{1}) \cup X(\beta_{2})$ fills $\Sigma$, that is every curve on $\Sigma$ has non-zero intersection with some curve from this set. It follows that $f$ must fix every curve on $\Sigma$. Therefore, $f \in Ker(\Sigma)$ and by assumption this is absurd. $\diamond$\\

\setlength{\parindent}{2em}

The link of $\alpha$ is either chain-connected, so that for any two of its vertices, $\beta_{1}$ and $\beta_{2}$, there is a sequence of curves $\beta_{1} = \delta_{1}, \delta_{2}, \ldots, \delta_{n} = \beta_{2}$ each distinct and disjoint from $\alpha$ and such that consecutive curves $\delta_{i}, \delta_{i+1}$ have minimal intersection, or is connected. By applying Lemma 16 inductively, we conclude that any two edges ending on $\alpha$ are prescribed the same automorphism of $\mathcal{C}(\Sigma)$ and that any such automorphism is induced by a mapping class. Since $\mathcal{C}(\Sigma)$ is connected, it follows that every edge is allocated the \tit{same} such automorphism $\Phi$.

All we need do now is verify that this automorphism is equal to $\phi$ everywhere. To do this, we only need to remark that any curve $\alpha$ spans an edge with a second curve $\beta$. This edge is prescribed the automorphism $\Phi$ which, by construction, agrees with $\phi$ on both $X(\alpha)$ and $X(\beta)$. In particular, $\Phi$ agrees with $\phi$ on $X(\beta)$ which contains $\alpha$. This completes one proof of Theorem 1.\\\\

\setlength{\parindent}{0em}
\setlength{\parskip}{0.5em plus 0.5em minus 0.5em}

\tbf{References.\\}

[BehrMar] J. Behrstock, D. Margalit, \tit{Curve complexes and finite index subgroups of mapping class groups} : to appear in Geometriae Dedicata.

[BelMar1] R. W. Bell, D. Margalit, \tit{Braid groups and the co-Hopfian property} : to appear in the Journal of Algebra.

[BelMar2] R. W. Bell, D. Margalit, \tit{Injections of Artin groups} : Preprint, arXiv: math.GR/0501051.

[Birm] J. Birman, \tit{Braids, links and mapping class groups} : Princeton University Press, Princeton, N. J., 1974, Annals of Mathematical Studies \tbf{82}.

[BreMar] T. E. Brendle, D. Margalit, \tit{Commensurations of the Johnson kernel} : Geometry \& Topology \tbf{8} (2004) 1361--1384.

[BroCMin] J. F. Brock, R. D. Canary, Y. N. Minsky, \tit{The classification of Kleinian surface groups, II: The ending lamination conjecture} : Preprint, arXiv: 
math.GT/0412006.

[D] M. Dehn, \tit{Papers on group theory and topology} : ed. J. Stilwell, Springer (1987).

[FIva] B. Farb, N. V. Ivanov, \tit{The Torelli geometry and its applications} : Mathematical Research Letters \tbf{12} (2005) 293--301.

[H] W. J. Harvey, \tit{Boundary structure of the modular group} : in ``Riemann surfaces and related topics: Proceedings of the 1978 Stony Brook Conference'' (ed. I. Kra, B. Maskit), Annals of Mathematical Studies No. 97, Princeton University Press (1981) 245--251.

[HK] W. J. Harvey, M. Korkmaz, \tit{Homomorphisms from mapping class groups} : Bulletin of the London Mathematical Society \tbf{37} No. 2 (2005) 275--284.

[Irm1] E. Irmak, \tit{Superinjective simplicial maps of complexes of curves and injective
homomorphisms of subgroups of mapping class groups} : Topology \tbf{43} No. 3 (2004) 513--541.

[Irm2] E. Irmak, \tit{Superinjective simplicial maps of complexes of curves and injective homomorphisms of subgroups of mapping class groups II} : Topology and Its Applications \tbf{153} No. 8 (2006) 1309--1340.

[Irm3] E. Irmak, \tit{Complexes of non-separating curves and mapping class groups} : Preprint, arXiv:math.GT/0407285.

[Iva1] N. V. Ivanov, \tit{Automorphisms of complexes of curves and Teichm\"uller spaces} : International Mathematics Research Notices \tbf{14} (1997) 651--666.

[Iva2] N. V. Ivanov, \tit{Subgroups of Teichm\"uller modular groups} : Translations of Mathematical Monographs \tbf{115}, American Mathematical Society (1992).

[IvaMcCar] N. V. Ivanov, J. D. McCarthy, \tit{On injective homomorphisms between Teichm\"uller modular groups I} : Inventiones Mathematicae \tbf{135} (1999) 425--486.

[K] M. Korkmaz, \tit{Automorphisms of complexes of curves on punctured spheres and on punctured tori} : Topology and its Applications \tbf{95} (1999) 85--111.

[L] F. Luo, \tit{Automorphisms of the complex of curves} : Topology \tbf{39} (2000) 283--298.

[Mar] D. Margalit, \tit{Automorphisms of the pants complex} : Duke Mathematical Journal \tbf{121} No. 3 (2004) 457--479.

[Sch] P. Schmutz Schaller, \tit{Mapping class groups of hyperbolic surfaces and automorphism groups of graphs} : Compositio Mathematica \tbf{122} No. 3 (2000) 243--260.

[Sha] K. J. Shackleton, \tit{Aspects of the curve complex and the mapping class group} : University of Southampton PhD thesis (2005).

[V] O. Viro, \tit{Links, two-sheeted branched coverings and braids} : Soviet Math, Sbornik \tbf{87} (2) (1972) 216--228.

\end{document}